\documentclass[12pt]{article}
\usepackage[numbers]{natbib}
\usepackage{amssymb}
\begin{document}

\title{\bf The arrow of time, complexity and the scale free analysis}

\author{
Dhurjati Prasad Datta\thanks{email:dp${_-}$datta@yahoo.com} 
\hspace{.1cm}and Santanu Raut\\
Department of Mathematics, 
University of North Bengal \\
 P.O. Siliguri,
West Bengal, Pin: 734430, India } 
\date{}
\maketitle

\baselineskip =20pt

\begin{abstract}
The origin of complex structures, randomness, and irreversibility are analyzed in the 
scale free SL(2,R) analysis, which is an extension of the ordinary analysis based on the 
recently uncovered scale free $C^{2^n-1}$ solutions to linear ordinary differential 
equations. The role of an intelligent decision making is discussed. We offer an 
explanation of the recently observed universal renormalization group dynamics at the 
edge of chaos in logistic maps. The present formalism is also applied to give a first 
principle explanation of 1/$f$ noise in electrical circuits and solid state devices. Its 
relevance to heavy tailed (hyperbolic) distributions is pointed out.  
\end{abstract}

\begin{center}
\bf Chaos, Solitons and Fractals 28, 581-589, (2006)
\end{center}
\newpage

\section{Introduction}

\par The origin of the arrow of time is still considered to be a puzzling problem in 
theoretical physics \cite{zeh}. Another difficult problem is the 1/$f$ noise 
\cite{press}, a footprint of complexity. The ubiquity of 1/$f$-like noise in diverse 
natural and biological processes seems, in particular, to signal to certain key, but 
still not clearly understood, dynamical principles that might be at work universally at 
the heart of any given dynamical process. Recently, there seems to have been an emerging 
urge in literature \cite{west} for a new principle for understanding complex, 
intrinsically irreversible, processes in nature.
\par In this paper, we argue that  the recently discovered finitely differentiable {\em 
scale free} solutions to the simplest scale free ordinary differential equation (ODE) 
\cite{dp1}-\cite{dp4}
\begin{equation}\label{ss}
t{\frac{{\rm d}\tau}{{\rm d}t}}=\tau
\end{equation}

\noindent should be able to offer an ideal framework, analogous to one advocated in 
Ref.\cite{west}, for complex phenomena. In particular, we reemphasize  how both the 
above problems might acquire a natural explanation in the class of such solutions. Both 
the problems were originally discussed in \cite{dp1,dp2,dp3}. The present discussion 
based on the {\em new exact} class of {\em time asymmetric} $C^{2^n-1}$ solutions to 
eq(\ref{ss})  \cite{dp4} is, however, expected to put the previous analyses in a 
stronger foundation. The simplest ode (\ref{ss}) will also be shown to enjoy a host of 
{\em new} dynamical features shedding  new lights in the definitions of living and 
nonliving systems. The relevance of the present formalism to complex systems has already 
been pointed out in literature \cite{weiss}.

\par To summarize briefly the previous results,    we present, in Ref\cite{dp1}, a novel 
{\em dynamical} treatment of linear ODEs when the time (i.e. the independent real) 
variable $t$ is {\em assumed} to have a random element. We show how a judicious use of 
the golden mean partition of unity, $\nu^2 + \nu = 1,\,\nu = (\sqrt 5-1)/2$, not only 
allows time to undergo random changes (flips) by inversions, $t_-\rightarrow t_-^{-1} = 
t_+,\,t_{\pm}=1\pm \eta$, in the vicinity of an instant $t=1$ (say), but also unveils, 
for the first time,  the possible existence of a class of random, second derivative 
discontinuous, scale free solutions to eq(\ref{ss}). One of the major aims of \cite{dp1} 
is to explore the especially unique role of the golden mean number in this extended 
framework of calculus accommodating inversions as a valid mode of changes (increments) 
besides ordinary translations. The freedom of random inversions provides a dynamic, 
evolutionary character to the second derivative discontinuous solutions, with a 
privileged sense of time's arrow, represented graphically by the slow, progressive 
unfolding of the approximants of the golden mean continued fraction. We also give a 
derivation of a universal probability distribution for these intrinsic random 
fluctuations. These are further clarified in \cite{dp2} where we treated, in a sense, 
the `inverse problem'. We show that the existence of a nontrivial solution of the form 
$\tau(t)=t (1+ \phi(t)),\,\phi(t)=t_1\tau(t_1^{-1})\neq$ an ordinary constant, 
$t_1=\epsilon t,\, \epsilon$ being a scaling parameter, to eq(\ref{ss}) could {\em only} 
be interpreted in an extended framework of calculus accommodating random inversions as a 
valid mode of changes for a {\em genuine} set of infinitesimals \cite{robin}. We also 
point out how this extended {\em dynamical calculus} provides a natural framework for 
the explanation of the origin of the generic 1/$f$ spectrum. These ideas become more 
precise and concrete in \cite{dp3} where we present an explicit solution of the form 
$\tau(t_-)=1/t_+$ to eq(\ref{ss}). This solution, though approximate ($\sim$O($\eta^2$) 
=0) , in the ordinary real number system $R$, is, however, {\em exact} in an nonstandard 
real number set $\bf R$. We also discuss an application of these {\em dynamical real 
numbers} to cell division and present an interesting {\em intelligent} version of the 
Newton' first law of motion. Finally, in \cite{dp4}, we show that the `approximate' 
solution of \cite{dp3} is in fact {\em generic}, in the sense that  the more accurate, 
in fact the exact solution, derived by generating successive self- similar corrections 
to an initially approximate solution, fails to yield  the exact  solution even in the 
limit of infinite number of iterations, thereby proving in a rigorous manner, our 
original contention that {\em any } measurement, even in the classical sense, {\em is 
intrinsically} uncertain. 
\par The new results presented here are the following. Besides restating the above 
novels features (viz., randomness and time irreversibility) in the context of the new 
`exact' but nevertheless intrinsically uncertain solution of \cite{dp4} in a much more 
greater detail, we examine the interesting problem of the origin of power law 
sensitivity in the logistic map at the edge of chaos \cite{tsallis2}, the 1/$f$ noise of 
electrical circuites and semiconductor devices \cite{voss}, and the long tailed 
hyperbolic type distributions in this extended framework of dynamics. 
 \par The paper is organized as follows. We first review the derivation of exact, higher 
derivative discontinuous ($C^{2^n-1}$) solutions of eq(\ref{ss}) in Sec.2. The new 
solution breaks the reflection symmetry ($t\rightarrow -t$) of the underlying ODE 
spontaneously \cite{dp3}. We also show here that besides these finitely differentiable 
($C^{2^n-1}$) time asymmetric solutions as well as the infinitely differentiable, time 
(reflection) symmetric standard solution, eq(\ref{ss}) possesses another {\em new} class 
of fluctuating solutions which are both infinitely differentiable and time symmetric. 
Because of these nontrivial classes of finitely and infinitely differentiable 
fluctuating solutions, a real variable $t$ can {\em undergo changes} not only by linear 
translations, but by inversions ($t \rightarrow 1/t$), in the neighbourhood of each real 
$t$. We next discuss how this defines a nonstandard extension of the real number system 
\cite{robin} in Sec.3. This also clarifies the origin of an {\em intrinsic randomness} 
at {\em as  fundamental a level as the real number system}.   The solutions are {\em 
intrinsically} stochastic  with a global ( irreversible) arrow of time (Sec.4).  
Consequently, every real number is identified with an equivalence class of a continuum 
of new, infinitesimally separated elements, which are in a state of random fluctuations. 
In Sec.5, we reinterpret,in the present framework, some new results on unimodal logistic 
map \cite{tsallis2,tsallis3} at the chaos threshold. In particular, we give a derivation 
of the $q$ exponential power law dynamics of the sensitivity to initial conditions. As 
another application of the late time variability of any dynamical system, arising from 
the ``primordial" (infinitesimal) {\em numerical fluctuations}, with generic $1/f$ 
spectrum, we discuss 1/$f$ noise in electrical circuits \cite{voss} in Sec.6.  In Sec.7, 
we show how a hyperbolic type distribution arises naturally at the asymptotic late time 
($t \rightarrow \infty$) limit even from a normally distributed variate.  

\section{Mathematical results}

\par Because of the novelty of the results, it is instructive to give a fairly complete 
derivation of such solutions \cite{dp4}. To this end, let us first construct the 
solution in the neighbourhood of $t=1$.  We need to introduce following notations.

\par Let $t_{n\pm}=1\pm \eta_n,\, t_0 \equiv t,\, 0<\eta_n<<1, \, 
\alpha_n=1+\epsilon_n,\, n=1,2,\ldots,\, {\rm and}\,\, 0 < \epsilon_n <1 $, such that 
$\epsilon_n \rightarrow 0, \, {\rm as} \, n \rightarrow \infty$. Next, we write 
$t_{n\pm}^\prime=1 \pm \alpha_n \eta_n^\prime$, so that $\alpha_n t_{n-} 
=t_{n-}^\prime$. Consequently, $\eta_n^\prime=\eta_n - {\frac{\epsilon_n}{\alpha_n}}$. 
Here, $\alpha_n$ (and $\epsilon_n$) are arbitrary scaling parameters. As will become 
clear,  natural parametrization is, however, $\epsilon _n=\epsilon _1^n$ \cite{dp4}. 
Further, $\eta_{n+1}=\alpha _n^2\eta^{\prime \,2}_n$.
\par To construct a nontrivial solution (with the initial condition $\tau(1)=1$), we 
begin with an initial approximate solution in the small scale variable $\eta_0$, and 
then obtain recursively self-similar corrections over smaller and smaller scale 
variables $\eta_0^2,\, \eta_0^4,\ldots$. To this end, let  

\begin{equation}\label{ns}
\tau (t)=
\cases{
\tau_- & {\rm if} $t\lessapprox 1$\cr
\tau_+ & {\rm if} $t\gtrapprox 1$},\,\,
\tau_{-}(t_-)=(1/t_+)f_{1-}(\eta_0),\, \tau_+(t_+)=t_+
\end{equation}

\noindent be an exact solution of eq(\ref{ss}). This is obviously true for the right 
hand component $\tau_+$.  Now the same will be true for the nontrivial component 
$\tau_-$ provided the correction factor $f_{1-}$ solves exactly the self-similar 
equation

\begin{equation}\label{ss1}
t_{1-}{\frac{{\rm d}f_{1-}}{{\rm d} t_{1-}}}=f_{1-} 
\end{equation}

\noindent in the smaller logarithmic variable $\ln t_{1-}^{-1}$, where 
$t_{1-}=1-\alpha_0^2 \eta_0^{\prime \, 2}\equiv 1-\eta_1$.  The exact (nontrivial part 
of the ) solution could thus be written recursively as 

\begin{equation}\label{recu}
\tau_-=C{\frac{1}{t_+}}{\frac{1}{t_{1+}^\prime}}\ldots
{\frac{1}{t_{{(n-1)}+}^{\prime}}}f_{n-}(\eta_n^\prime)
\end{equation}

\noindent where $f_{n}$ satisfies the $n$th generation self-similar equation  

\begin{equation}\label{ssn}
t_{n-}{\frac{{\rm d}f_{n-}}{{\rm d} t_{n-}}}=f_{n-}
\end{equation}

\noindent and $t_{n-}=1-\alpha_{n-1}^2 \eta_{n-1}^{\prime \, 2}\equiv 1-\eta_n$ so that 
$\eta_{n}=\alpha^2_{n-1}\eta^{\prime 2}_{n-1}$. We also note that $t_+^\prime= t_+$, 
since $\alpha_0=1$.

\par Plugging in the initial condition $\tau_{\pm} =1$ at $t_{\pm} =1$ (viz., $\eta_0 
=0$), one obtains finally the desired solution as 

\begin{equation}\label{ns1}
\tau_-=C{\frac{1}{t_+}}{\frac{1}{t_{1+}^\prime}}{\frac{1}{t_{2+}^{\prime}}}
\ldots, \, \, \, \tau_+=t_+
\end{equation}

\noindent where $C=t_{1+}^\prime(0)t_{2+}^\prime(0)\ldots$. Notice that $C\neq 1$, since 
$ \eta_1^\prime= -\epsilon_1/\alpha_1, \,  \eta_2^\prime=\epsilon_1^2 
-\epsilon_2/\alpha_2, $ etc, when $\eta_0 =0$. Notice also that it is easy to show that 
$\{f_{n-}\}$ is convergent and that $f_\infty = {\frac{{\rm d}f_{n-}}{{\rm d} 
t_{n-}}}|_\infty ={\frac{{\rm d}\tau}{{\rm d} t}}|_{t=1} = \tau(1)=1$, in the asymptotic 
limit $n \rightarrow \infty$ \cite{dp4}.

\par The salient features of this solution are the following.

1. The solution has discontinuous second derivative at $t=1$.  The said discontinuity is 
an effect of an {\em infinity} of nonzero  rescaling parameters $\epsilon_n$. For a 
finite set of $\epsilon_n$ ( or in the special case when $\epsilon_n=0,\, \forall n$), 
one gets back the standard solution. Moreover, the scale invariance is realized only in 
a one sided manner. The scaling $\alpha_n t_-=t_-^\prime$ does not mean $\alpha_n 
t_+=t_+^\prime$.

2. It also follows that the solution (\ref{ns1}) is indeed an exact solution of 
eq(\ref{ss}) when the ordinary real variable $t=1-\eta_0$ (near $t=1$) is replaced by 
the {\em fat} real variable ${\bf t_-^{-1}}=\Pi^\infty_0 t_{n+}^\prime$. The fat 
variable $\bf t$ leaves in $\bf R$, a nonstandard extension \cite{robin} of the ordinary 
real number set $R$, inhabiting {\em infinitesimal} scales (variables) $\ln t_n^\prime 
\approx \alpha_n \eta_n^\prime$. All these variables can be treated as independent 
because of the arbitrary scaling parameters $\epsilon_n$.

3. The scaling invariance of eq(\ref{ss}) tells also that, $t=1$ could be realized as $t 
\rightarrow t/t_0=1$, so that the nontrivial solution (\ref{ns1}) actually holds in the 
neighbourhood of every real number $t_0$, the 2nd derivative being discontinuous at 
$t=t_0$.  Combining the standard and the new solutions together, one can thus write down 
a more general class of solutions in the form 

\begin{equation}\label{gen}
\tau_g(t)=t(1 + \phi(t)), \, \, \phi(t)=\epsilon t^{-1} \tau(t_1),\,t_1=t/\epsilon 
\end{equation} 

\noindent where $\epsilon$ is another arbitrarily small scaling parameter. Note that  

\begin{equation}\label{gc}
t{\frac{{\rm d}\phi}{{\rm d}t}}=0
\end{equation}

\noindent because $\tau$ is an exact solution of eq(\ref{ss}) (one may  fix $\epsilon_n 
= \epsilon^{2^n}$, in view of remark 2). The 2nd derivative discontinuity of $\tau$, 
however, tells that $\phi$ can not be considered simply as an ordinary 
constant~\cite{dp1, dp2} (for an application see Sec.5).

4. Finally, it is easy to verify that eq(\ref{ss}) possesses yet another (nontrivial ) 
class of infinitely differentiable solutions of the form 

\begin{equation}\label{ns2}
\tau^\prime (t)=
\cases{
\tau_-^\prime & {\rm if} $t\lessapprox 1$\cr
\tau_+^\prime & {\rm if} $t\gtrapprox 1$},\,\,
\tau_{-}^\prime(t_-)=(1/t_+)f(\eta_0),\, \tau_+(t_+)^\prime=(1/t_-)f(\eta_0),
f(\eta_0)={\frac{1}{t_{1+}^\prime}}{\frac{1}{t_{2+}^{\prime}}}\ldots
\end{equation}
 
 \noindent which is, however, distinct from the standard solution. This extends the 
purely fluctuating time symmetric solution studied in ref\cite{dp3}. Note that the 
infinite differentiability is restored because of identical self similar corrections in 
$\tau_\pm^\prime$. However, as it should be evident from the above derivations, the 
iteration schemes for both $\tau_-^\prime$ and $\tau_+^\prime$ could be run 
independently with different sets of scaling factors $\epsilon_n$ and $\epsilon_n^ 
\prime$ respectively, leading again to second derivative discontinuity. Besides these 
second derivative discontinuous solutions,  eq(\ref{ss}), as shown in \cite{dp4}, also 
accommodates a larger class of $C^{2^n- 1}$ solutions. Consequently, the simple ode 
(\ref{ss}) accommodates indeed {\em an astonishingly rich set of solutions belonging to 
different differentiability classes}.

\section{Randomness}

\par Is randomness a fundamental principle (law) controlling our life and all natural 
processes? Or, is it simply a projection of our limitations in comprehending such 
complex phenomena? Is there any well defined boundary separating simple and complex?  
These and similar related questions on the actual status of randomness are being 
vigorously investigated in literature \cite{west, prigo}. The new class of discontinuous 
solutions sheds altogether a new light on the ontological status of randomness. (It is 
believed that randomness in quantum mechanics arises at a fundamental level. However, 
the Schr{\~o}dinger equation, the governing equation of any quantal state, is purely 
deterministic and time symmetric. The random behaviour is ascribed only through an 
extraneous hypothesis of a 'collapsed state' at the level of measurement (see, \cite 
{west}))

\par Let us note that eq(\ref{ss}) is the simplest ODE, $t$ being an ordinary real 
variable. In the framework of the conventional analysis, one can not, in any way, 
expect, at such an elementary level, a random behaviour in its, so called unique 
(Picard's), solution. (The chaos and unpredictability could arise only in the presence 
of explicit nonlinearity in higher order ODEs.) However, the nontrivial scaling, 
exploited in Sec.2, along with the initial ansatz (\ref{ns}), reveals not only the self 
similarity of $C^{2^n-1}$ solutions over scales $\eta^{2^n}$, but also exposes a subtle 
role of {\em decision making} and randomness in generating nontrivial late $t$ behaviour 
of the solution.    

\par A basic assumption in the framework of the standard calculus is that a real 
variable $t$ changes by linear translation only. Further, $t$ assumes (attains) every 
real number exactly. However, in every computational problem within a well specified 
error bar, a real number is determined only up to a finite degree of accuracy 
$\epsilon_0$ ,say. Suppose,  for example, in a computation, a real variable $t$ is 
determined upto an accuracy of $\pm 0.01$, so that $t$  here effectively   stands for 
the set $t_\epsilon\equiv\{t \pm \epsilon,\,\epsilon < \epsilon_0 = 0.01\}$ with 
cardinality $c$ of the continuum. We call $\epsilon$ an `infinitesimally small' real 
number (variable). Now, any laboratory computational problem (experiment) is run only 
over a finite time span, and the influence of such infinitesimally small $\epsilon$'s, 
being insignificantly small, could in fact be disregarded. Consequently, the variable 
$t$ could be written near 1, for instance, as $t_1=1+\eta$, where $\eta$ is an ordinary 
real variable close to 0, having 'exact' values as long as one disregards infinitesimal 
numbers ($< \epsilon_0$) due to {\em practical limitations}. At the level of 
mathematical analysis, such practical limitations being indicative only of natural 
(physical/biological) imperfections should not jeopardy the existence of an abstract 
theory of sets, real number system, calculus and so on, shaping the logical framework 
for an exact and deterministic understanding of natural processes. That such an attempt 
would remain as an unfulfilled dream is not surprising in view of the new class of 
discontinuous solutions. As shown in detail in ref\cite{dp4}( see also 
\cite{dp1,dp2,dp3}) the $C^{2^n-1}$ solutions shows that the real number set $R$ should 
actually be identified with an nonstandard real number set ${\bf R}$ \cite{robin} so 
that every real number $t$  {\em is a fat} real  (hyperreal ) number ${\bf t}$, which 
means that $t \equiv \bf t$. Accordingly, a real number $t$ could not be represented 
simply by a structurless point, but in fact is embedded in a sea of irreducible 
fluctuations of infinitesimally small numbers ${\bf t}=tt_f$, $t_f= 1+ \phi$ denoting 
{\em universal} random corrections from infinitesimals $\phi$ (c.f., eqn(\ref{gen})). 
The origin of randomness is obviously tied to the {\em freedom} of injecting an infinite 
sequence of arbitrary scaling parameters $\epsilon_n $ into the $C^{2^n-1}$ solutions  
because of scale invariance of eq(\ref{ss}), introducing small scale uncertainty 
(indeterminacy)  in the original variable. (Note that scaling at each stage introduces a 
degree of uncertainty in the original variable viz., $\delta \eta_1 = |\eta_1^\prime - 
\eta_1|=\epsilon_1/\alpha_1$ and so on.) It is also shown that the set of infinitesimals 
$\bf 0=\{\pm \phi\}$ has a Cantor set like structure viz., discrete, dust like points 
separated by voids of all possible sizes-- and having the cardinality $2^c$ 
\cite{dp1,dp2}. Accordingly, an infinitesimal variable could change (within an 
infinitesimal neighbourhood of a point, say, 1) only by discrete {\em jumps} (inversions 
) of the  form $t_-=t_+^{-\alpha},\, t_\pm= 1\pm\phi$, to cross the gulf of emptiness,  
length of jumps being arbitrary because of an arbitrary  $\alpha$.      
Note that, the value of a small real number $\eta_0$ is uncertain not only upto 
O($\eta_0^2$), but also because of the arbitrary parameter $\epsilon$. Note also that 
the solution (\ref{ns1}) proves explicitly that inversion is also a valid mode of change 
for a real variable, at least in an {\em infinitesimal} neighbourhood of an ordinary 
real number. We note further that two solutions $\tau_g$ and $\tau_s$ are 
indistinguishable for $t \sim$ O(1) and $\eta_0^2 <<1$. However, for a sufficiently 
large $t \sim$ O($\epsilon^{-1})(\equiv$ O($\eta_0^{-2}$)), the behaviours of two 
solutions would clearly be different. Finally, the order of discontinuity could be 
``controlled" by an application of an {\em intelligent decision} invoking a nonzero 
$\epsilon _n \neq 0,\,\epsilon _m =0,\,m<n$ only at the $n$th level of iteration. This 
{\em freedom of decision making} could either be utilized at a pretty early stage of 
iterations, for instance, $n=1$, say, making the system corresponding to the ode 
(\ref{ss}) {\em fully intelligent}, or be {\em postponed } indefinitely ($n= \infty$) 
reproducing the standard Picard's solution for a material ( non-intelligent) system 
\cite{dp3}. We note that this randomness and potential  intelligence being intrinsic 
properties of the solutions of eq(\ref{ss}) {\em are indelibly   rooted to the real 
number system and hence could not simply be interpreted as due to some coarse graining 
effect} analogous to thermodynamics and statistical mechanics. Randomness at a 
fundamental level arises also  in the E- infinity theory \cite{nash1} ( for other 
approaches, see \cite{west}).

\par Before closing this section, let us restate the  point of view that is being 
emerged out of our analysis. Because of historical reasons, current thought processes in 
any scientific fields of research--- physical, chemical, biological, financial, 
sociological,  ... is largely shaped by the role model- Physics(!), which conceives 
nature as composed purely of inanimate matter and force fields. Consequently, life 
together with its complex ramifications get reduced to the level of derived concepts - 
coarse-grained, in some sense, from `more fundamental'(?) material principles (recall 
the present state of art of the main-stream biological (biophysical) research!). Even in 
a purely (lifeless) physical world the origin of randomness from some sort of 
coarse-graining presupposes our {\em practical} limitations in exact measurability (in 
classical sense) (as in thermodynamics and statistical mechanics, see above). Hence both 
life-like properties and randomness seem to have been relegated to the level of 
secondary derived concepts, arisen out purely of inanimate material principles and /or 
``our ignorance" \cite{prigo} in comprehending the world, rather than being linked 
directly to some intrinsic dynamical properties (aspects and/or quality) of  life (and 
nature) itself \cite{west}. This state of affairs does not yet have undergone any 
dramatic  changes with the introduction, as remarked already, of quantum mechanics and 
quantum field theories, because, the status of randomness in a quantal theory is not yet 
clear. The new mathematical results  explored here ( and in \cite{dp1}-\cite{dp4}) would 
likely to have some profound implications in this regard. With intelligence (and 
decision making) emerging as a fundamentally new ingredient (degree of freedom) from the 
mathematical analysis, the traditional framework of a physical theory viz., space, time, 
matter, energy ( or in a relativistic theory, spacetime and energy) might, in future, be 
extended and replaced by a {\em truly dynamical} framework consisting of {\em 
intelligence, space, time, matter, energy} as envisioned in Ref.\cite{west}. To explore 
the dynamic properties of the new solutions further, we now examine the origin of time's 
arrow  in the following.

\section{Reflection symmetry breaking and time}
\par It is well known that {\em time} is {\em directed}, that is to say, we all have a 
sense of a forward moving time. The problem of time asymmetry [1] points to a 
fundamental dichotomy between the (Newtonian) `time' in physics and mathematics and that 
of our (objective) experiences. The Newtonian time is non-directed. There is no way to 
distinguish between a space like variable $x$ with a time variable $t$. Further, 
all the fundamental equations of physics are time reversal symmetric. However, the 
existence of $C^{2^n - 1}$solutions of eq(\ref{ss}) presents us with a new scenario! 
One is now obliged to re-examine the conventionally accepted standard notions under this 
new light. As mentioned already, we show here that time {\em does} indeed has an arrow, 
which is inherited,  not only by all (physical /biological /social) dynamical systems, 
but it is also indelible inscribed even to a real number. The concept of time 
thus turns out to be more fundamental compared to space and may even be considered at 
par with the real number system (and hence to the existence of intelligence as a 
fundamental entity)!

\par To see how a time sense is attached to $C^{2^n}$ solutions, we recall first that 
the variable $t$ is, in general, non-dynamical, and need not denote the (forward) flow 
of time. In fact, it simply behaves as a labeling parameter. Further, the inversion 
$t_-\leftrightarrow t_-^{-1}= t_+$ may at most be considered as a reversible random 
fluctuation between $t_\pm$ (for a given $\eta>0)$ with equal probability 1/2. In the 
usual treatment of ordinary calculus and classical dynamics, $t$ is a non-random 
ordinary variable, and the above inversion reduces to the symmetry of eq(\ref{ss}) under 
reversal of sign (parity) $t\rightarrow -t$. The infinitely differentiable standard  
solution $\tau_s(t)$, written (in the notation of eq(\ref{ns})) as $\tau_{s-}=t_-, \, 
\tau_{s+}=t_+ ,\,(\tau_s(1)=1)$ is obviously symmetric under this reversible inversion 
($\eta\rightarrow -\eta$). The nontrivial solution $\tau(t)$ in eq(\ref{ns1}), however, 
constitutes an explicit example where {\em this parity invariance is dynamically 
broken}, viz.; when the inversion is realized in an irreversible (one-sided, directed) 
sense. 
\par To state the above more precisely, let $P: P t_\pm =t_\mp$ denote the reflection 
transformation near $t=1$ ($P\eta=-\eta$ near $\eta=0$). Clearly, eq(\ref{ss}) is parity 
symmetric. So is the standard solution $\tau_{s\pm} = t_\pm$ (since $P \tau_s = 
\tau_s$). However, the solution (\ref{ns1}) breaks this discrete symmetry spontaneously: 
$\tau_-^P=P \tau_+ = t_-,\, \tau_+^P = P \tau_- = C{\frac{1}{t_-}} 
{\frac{1}{t_{1+}^\prime}} {\frac{1}{t_{2+}^{\prime}}}\ldots$, which is of course a 
solution of eq(\ref{ss}), but clearly differs from the original solution, $\tau_\pm^P 
\neq \tau_\pm$.   
\par To see more clearly how this one-sided inversion is realized, let $t\rightarrow 
1^-$ from the initial point $t\approx 0$. Then at {\em a} point in the infinitesimal 
neighbourhood of  $t_-\lessapprox 1$, the solution $\tau_-$ carries (transfers) $t_-$ 
instantaneously to $t_+$ by an inversion $\tau_-\approx 1/t_+$ (we disregard here the 
O$(\eta_0^2)$ and lower order self similar fluctuations) and subsequently the solution 
follows the (standard) path $\tau_+=t_+$ in the small scale variable $\eta\gtrapprox 0$, 
as it is now free to follow the standard path  till it grows to O($\lessapprox 1$), when 
second order transition to the next smaller scale variable by inversion becomes 
permissible: $ \eta=1/\eta_+, \eta_+=1+\bar \eta,\,\bar\eta\gtrapprox 0$, and so on. 
Clearly, the generic pattern of (irreversible) (time asymmetric) evolution in $\tau(t)$ 
over smaller and smaller scales resembles more and more closely the infinite continued 
fraction of the golden mean: $\tau_-(t)=t, \,0<t\uparrow \lessapprox 1,\, \rightarrow 
\tau_-=1/(1+\eta), \, 0<\eta\uparrow \lessapprox 1,\, \rightarrow \tau_-=1/(1+1/(1+\bar 
\eta)),\,\bar\eta\gtrapprox 0$, and so, $\tau_-(t)\rightarrow \nu,\, \nu = (\sqrt 
5-1)/2$, the golden mean, as $t \rightarrow \infty$. Here, $t \uparrow $ means that $t$ 
is an increasing variable. We note that the (macroscopic) variable $t$ is {\em 
reversible} as long as $t \sim 1$. Subsequently, this parity symmetry is broken by a 
random inversion, leading the evolution irreversibly to a smaller scale $0<\eta<<1$, so 
that $\eta$ has only one way to change (viz., $\eta\rightarrow 1^-$), since the 
transition $t_+\rightarrow t_-$ ( the parity reversed solution) is improbable (Note that 
the growth of $\eta_0$ in this case is obstructed due to divergence of $1/t_-$ as 
$\eta_0 \rightarrow 1^-$). We also note that the solution which is written only in the 
neighbourhood of $t=1$, could be extended to the above more general pattern of a `time 
irreversible' solution, accommodating an infinite (discrete) set of points where the 
higher order smoothness of this general function meets obstructions in the form of 2nd 
(or higher) derivative discontinuity, forcing it to change stochastically its path down 
the cascades of the infinite continued fraction of the golden mean. Consequently,  {\bf 
time} is indeed {\em rediscovered} as one with its true elements viz., {\em stochastic 
irreversibility }, rather than being considered simply as a labeling parameter. 
Analogous ideas are also advocated by Prigogine \cite{prigo} and El Naschie \cite{nash1, 
nash2}.       

\par The purely fluctuating solution eq(\ref{ns2}) is, however, reflection (time) 
symmetric, as could be verified easily. However, because of intrinsic randomness, this 
solution would also fail to return an initial value exactly, viz., $t_- \rightarrow t_+ 
\rightarrow t_+^\prime, \, t_+^\prime$ being only approximately equal to $t_-$. This 
slight mismatch between the initial and final values would definitely induce a symmetry 
breaking, thereby realizing the time asymmetric solution in due course. Consequently, 
{\em a variable $\tau$ is intrinsically time-like if and only if it corresponds to the 
solution (\ref{ns1}) which breaks the reflection symmetry spontaneously.} Obviously, 
this definition also applies to the fluctuating solution (\ref{ns2}) with different sets 
of scaling parameters breaking again the reflection symmetry of the underlying equation.  

\section{Logistic Maps}
\par Let us now point out an interesting application of the above formalism. There have 
been some new, previously unexposed, asymptotic scaling properties of the iterates of 
unimodal logistic maps at the edge of chaos \cite{tsallis2, tsallis3}.   For 
definiteness, we consider here only the dynamics of the sensitivity $\xi_t$ to the 
initial conditions for large iteration time $t$, at the chaos threshold $\mu =  
\mu_\infty = 1.40115...$ of the map $f_\mu(x) = 1-\mu |x|^2, \, -1 \leq x \leq 1$ (for 
notations see \cite{tsallis2}). For a sufficiently large $t$, the ordinary exponential 
behaviour of sensitivity gives away to a power law behaviour, and is shown to have the  
$q-$exponential form given by

\begin{equation}\label{sens1}
\xi_t= {\rm exp}_q(\lambda_q t) \equiv [1 +(1-q)\lambda_q t)]^{1/(1-q)}
\end{equation}

 \noindent The standard exponential dependence $\xi_t= {\rm exp}(\lambda _1 t)$ is 
retrieved at the limit $q \rightarrow 1$. The system is said to be strongly insensitive 
(sensitive) to initial conditions if $\lambda _1 <0$ ($\lambda _1>0$). The behaviour, 
however, gets altered at the edge of chaos. Using the Feigenbaum's RG doubling 
transformation $\hat R f(x)=\alpha f(f(x/\alpha))$ $n$ times to the fixed point map 
$g(x)$ viz., $g(x)=\hat R^n g(x) \equiv \alpha^n g^{2^n}(x/\alpha^n)$,  the values of 
the $q-$Lyapunov  coefficient $\lambda _q$ and $q$ are determined to be $\lambda _q = 
\ln \alpha/\ln 2 \,(>0),\, q= 1-\ln 2/\ln \alpha \,(<1),\, \alpha$ being one of the 
universal Feigenbaum  constants, $\alpha = 2.50290...$.  Consequently, the critical 
dynamics corresponds to {\em weak chaos}. Notice that the dynamics at the chaos 
threshold, being the most prominent and readily accessible to numerical experiments,  
among the critical points of a quadratic map, reveals a {\em universal} concerted 
behaviour, described by the fixed point solution of the RG doubling transformation. 
One, however, needs to explain the origin of weak chaos, i.e., the system being weakly 
sensitive, rather than being only weakly insensitive, to initial conditions, at the 
chaos threshold, when the system approaches the threshold from the left through period 
doubling route ( for $\mu < \mu _\infty$).   
 
\par To reinterpret above observations in the present context, we note first of all that 
the fixed point equation $\tilde g(x)= g(x)$ where $\tilde g(x) = \alpha 
g(g(x/\alpha))$, is a solution of eq(\ref{ss}), viz., $d\tilde g/d g=\tilde g/g$. 
Consequently, one expects that the critical dynamics of unimodal quadratic maps at the 
onset of chaos would be linked directly  to our nontrivial solutions to eq(\ref{ss}). As 
an explicit example, let us investigate the origin of the q-exponential type power law 
dynamics from the ordinary exponential one 

\begin{equation}\label{sens2}
{\frac{{\rm d}\xi}{{\rm d}t}}=\lambda_1 \xi
\end{equation}

\noindent As the critical point is approached from left (say), the Lyapunov exponent 
$\lambda_1 \rightarrow 0^-$, and hence gets replaced, in our extended framework, by a  
`dynamic' infinitesimal (ordinary zero is an equivalence class of infinitesimals) of the 
form $\lambda_1 \rightarrow -\epsilon \lambda_p \phi(t_1)$, where   $\epsilon \, (>0)$ 
is an infinitesimal ($\epsilon\neq 0,\,\epsilon^2=$O(0)) scaling parameter, $\phi = t_1 
\tau(t_1^{-1})$ ($t_1 =  \epsilon  t$), and $\lambda_p >0, \, p=p(\epsilon)$ are two 
generalized constants (c.f., eq(\ref{gc})), both of O(1).  However, when the variation 
in the intrinsic time like variable $\phi$ becomes relevant at the scale $t \sim 
1/\epsilon$,  that of $\lambda_p$ would be relevant only at a longer scale viz., $t \sim 
1/\epsilon^2$ or more and could be considered as an ordinary constant. The reason for 
introducing the (third) generalized constant (again with a slower variation than $\phi$ 
) $p$ will become clear below. At a critical point, Eq(\ref{sens2}) now reduces to the 
RG -like equation

\begin{equation}\label{sens3}
t_{1+}{\frac{{\rm d}\xi}{{\rm d}t_{1+}}}= \lambda_p \xi
\end{equation}
         
\noindent  To explain the derivation of the above equation let us proceed in steps. 

(i) The rescaled variable $t_1$, defined by $t_1 =\epsilon t$, is O(1) when $t \sim$ 
O($1/\epsilon)$. Consequently, the critical dynamics would be revealed only at {\em a 
sufficiently long time scale} $t_{1-}= 1-\eta \rightarrow 1^-$ i.e., $t \rightarrow 
(1/\epsilon)^-$.  That means {\em the dynamics at the chaos threshold needs to be probed 
in an extended framework, viz., in the sense of a limit as the control parameter $\mu $ 
approaches $\mu_\infty$ from left (say) through period-doubling cascade} instead of 
simply replacing $\mu$ by $\mu_\infty$ in the map and then iterating. In any 
computational problem, this extended framework is automatically realized, because of the 
inherent finite bit (decimal) representation of  a real number, such as $\mu_\infty$, 
exact value of which could only be approached recursively by increasing its accuracy. 
The infinitesimal $\epsilon$ (along with an O(1) variability  as represented by 
$\lambda_p \phi(t_1)$ ) then simply corresponds to the infinite trailing bits in the 
finitely represented real number e.g., $\mu_\infty$.   

(ii) The critical point equation (\ref{sens3}) follows when one makes use of the 
relation $d \ln t_{1-}= -d \ln t_{1+}$, which is valid for infinitesimal $\eta$ with 
O($\eta^2$)=0 \cite{dp3}. Clearly, the solution to the above equation is $\xi = [1+p 
\lambda _p \eta]^{1/p}$ which corresponds exactly to the $q$ exponential (\ref{sens1}) 
provided we choose $q=1-p,\, p \lambda_p =1$. 

(iii) The $q$ exponential solution is, however, valid not only in the `infinitesimal' 
neighbourhood of $t_1=1$, but for arbitrarily large $t_1$, because of the scale 
invariance of eq(\ref{sens3}). Indeed, writing $2^n t_1 = 1+\tilde t$, the $q$ 
exponential sensitivity takes the form $\xi = [1+(1-q) \lambda _q \tilde  t]^{1/(1-q)}, 
\, \tilde t$ being large. 

(iv) The scale factors $2^n$ correspond to the times to determine the trajectory 
positions $x_{2^n}$ of the logistic map with an initial position $x_{in-}$ 
\cite{tsallis2}. It follows therefore that the ratio of the sensitivities at times $2^n$ 
and $2^{n+1}$ viz., $\frac{\xi_{2^{n+1}}}{\xi_{2^n}} = (\frac{2^n t_1}{2^{n+1} 
t_1})^{1/p} = \alpha$ when $p$ remains constant at $1/p= \ln \alpha/\ln 2$ upto a time 
$T \approx 2^N,\, N \approx |\ln \xi_{in}/\ln \alpha|,\, \xi_{in}$ being the initial 
value of sensitivity, corresponding to the initial iterate $n =0$. 

(v) Indeed, to obtain the later estimate, we note that  $\xi_{in}^p = t_1$, which 
follows from  eq(\ref{sens3}). For a non-zero (sufficiently large) $n$, it now follows 
that $\xi_n \equiv (2^n \xi_{in}^p)^{1/p} =\alpha^n$  which translates to $ p = 
\frac{\ln 2}{\ln \alpha}(1 + \frac{\ln \xi_{in}}{n \ln\alpha})$. Consequently, a 
possible variation in $p$ would be revealed only when $n \approx |\ln \xi_{in}/\ln 
\alpha|$, as claimed. Note  that $q$ exponential form in the neighbourhood of $t_1 = 1$  
belongs to the class of solutions (\ref{ns2}) provided the generalized constant $p$ is 
given by $1/p =  \ln f/\ln t_1 - 1,\,(d p/\ln t_1) = 0$. 
 
\par The change in sensitivity from strongly insensitive case ($\lambda_1 <0$) to weakly 
sensitive ($\lambda_q>0,\,q<1$) power law behaviour is thus explained as an effect of 
nontrivial infinitesimals and associated inversion $t_{1-}t_{1+} \approx 1$ in the 
infinitesimal neighbourhood of $t_1=1$. The origin of Feigenbaum's constant $\alpha$ ( 
notice our use of $\xi_n =\alpha^n$ \cite{tsallis2,tsallis3} in the above derivation) in 
the present formalism along with other relevant issues will be considered elsewhere. 
Another interesting problem is to identify the golden mean number $\nu$ in the critical 
dynamics.    

\section{1/$f$ spectrum}

\par The relevance of higher derivative discontinuous solutions to the origin in $1/f$ 
noise problem have been discussed in detail in \cite{dp1,dp2}. We note here that in the 
extended framework of a dynamical theory, accommodating these solutions, any physical,  
variable $t$, say time, is replaced by $t^{1+\sigma},\,\sigma=\ln (1+\phi)/\ln t$ (c.f., 
Sec.3) where $\sigma$ typically is small O($\epsilon$) for any $\epsilon>0$. The nonzero 
exponent $\sigma$  introduces a small stochastic fluctuations over the ordinary (time) 
variable $t$. Clearly, these small scale stochastic fluctuations, existing purely in the 
real number system, would remain insignificant for any terrestrial (laboratory)  {\em 
inanimate} system  which persists over a moderate time scale, such as the motion of a 
(classical) particle under gravity. However, even for  simple  electrical circuits where 
the voltage fluctuation spectrum $S_V(f)$ is known to vary proportionally with the 
thermal fluctuation spectrum $S_T(f)$, the origin of 1/$f$ noise as observed in 
\cite{voss} could be naturally ascribed to the $C^{2^n-1}$ solutions of

\begin{equation}\label{th}
c{\frac{{\rm d}T}{{\rm d}t}}=-g(T-T_0)
\end{equation}

\noindent  This equation describes the macroscopic (equilibrium) variations of the 
temperature ($T$) of a resistive system with heat capacity $c$, coupled through a 
thermal conductance $g$ to a heat source  at temperature $T_0$ \cite{voss}. According to 
the conventional knowledge one does not expect 1/$f$ spectrum from such a simple, purely 
deterministic, linear uniscale system. One needs, in fact, to consider extraneous 
nonlinear effects from environment to explain the origin of the generic 1/$f$ 
fluctuations.     

\par However, according to the   present analysis, even this simple system would behave 
stochastically because of small scale, intrinsic fluctuations in the time variable $t$.  
These scale free fluctuations could influence the late time behaviour of the system 
provided {\em the system is 'allowed' to survive} over a period $t>> 
1/\epsilon,\,\epsilon = g/c$.  Notice that ordinarily a system following purely 
eq(\ref{th}) is assumed to  relax to the equilibrium temperature $T_0$ after $t \approx 
1/\epsilon$. The late time variability that is observed in any  resistive system is then 
ascribed normally to the complex nature of the resistive medium and/or (nonlinear) 
interactions with environment \cite{voss}, asking for an explicit modeling.   The 
generic observation of $1/f$ fluctuations in metal films and semiconductors still eludes 
a universal explanation for its microscopic origin in the framework of conventional 
dynamical theories \cite{planat1}. 

\par In view of $C^{2^n-1}$ solutions, we now have an extended framework to re-examine 
the above problem. The solution of eq(\ref{th}) now have the form $T(t)-T_0=t^\sigma 
e^{-\epsilon t}, \, \sigma $ being a small fluctuating variable. As noted already, this      
random exponent would lead to small scale stochastic modulations over the (mean) 
macroscopic decay mode, as observed in physical systems. These small scale (power law) 
fluctuations would persists even far beyond the ordinary relaxation time. Accordingly, a 
time series of temperature fluctuations ($T_f(t)=(T(t)-T_0)e^\tau$) recorded over a 
period of a few decades (1 to $10^4$, say), in the unit of the dimensionless time  $\tau 
= \epsilon t$, would reveal a scale free 1/$f$ type variability. For, the two point 
autocorrelation  function of this intrinsic fluctuations has the form $C(t)=< 
T_f(t)T_f(0)> =c<T_f(t)> =c<t^\sigma>\approx ct^{<\sigma>},\, T_f(0)=c$, where 
$<\sigma>$ is the expectation value of the random exponent and $c \sim$O(1) is the 
initial (background) noise in the system. The associated probability distribution would 
have a generic late $t$ behaviour, resembling infinitely divisible Levi type 
distributions \cite{levi} (see below). The corresponding power spectrum of this 
stochastic scale free fluctuation  is given by $S(f)\sim 1/f^{1-<\sigma>}$. We note that 
an arbitrarily small  nonzero $\sigma$ is sufficient to generate a $1/f$-like 
fluctuation. In other words, intrinsically random, infinitesimal scales associated with 
the time variable $t$ could act as a perennial  source of small scale fluctuations 
leading to the universal low frequency 1/$f$ spectrum. However, an accurate 
determination of the intensity of the fluctuations (viz, the constant of proportionality 
in the observed spectrum $S_V(f) \propto V^2/f^a,\,a \approx 1$) may require further 
work \cite{voss,planat1}. The relevance of number theory to 1/$f$ noise problem is also 
pointed out by Planat \cite{planat2}. El Naschie  suggested that the exponent $\beta$ of 
the $1/f^{\beta}$ noise to semiconductors would be related to the golden mean in the 
framework of the E-infinity theory \cite{nash3}. 

\section{Hyperbolic distribution}
 
\par In ref \cite{dp1} we show that the scale free infinitesimal fluctuations follow a 
nongaussian, Bramewell-Holdsworth-Pinton (BHP) \cite{bhp} distribution. Here we show how 
a hyperbolic, power law tail gets superposed generically in {\em any} distribution when 
the concerned random variate is assumed to leave in $\bf R$. To see this it suffices to 
consider only a normally distributed variate, because by the central limit theorem a 
normal probability density acts as the attractor for any probability density with finite 
moments. Let $t$ be a zero mean normal variate, with unit standard deviation. The 
corresponding fat variate could be written as ${\bf t}^2= t^2 + \phi,\,\phi= \epsilon(t) 
\ln t^2$, being a random infinitesimal satisfying eq(\ref{gc}) in the logarithmic 
variable $\ln t$ (c.f., Appendix). Consequently, the normal density function $\propto 
e^{-t^2/2}$ gets a generic power law tail $t^{-\epsilon} e^{-t^2/2}$. This generic power 
law tail in the present extended formalism should become important in the future studies 
on the statistics of rare events. We close with the remark that occasional detections of 
exceptional events in an experiment of, for instance, a normal variate over a prolonged 
period could be explained by this slowly varying tail. Note that the power law 
variability would become visible only in the asymptotic limits ($t \rightarrow 
\pm\infty$) because of an infinitesimal $\epsilon$, so that $\phi$ remains vanishingly 
small in any laboratory experiments over a finite time scale.       

\section*{Appendix}

A continuously differentiable function $f(t)$ of a real variable could be defined as an 
integral of the ODE ${\frac{{\rm d}x}{{\rm d}t}} = f^{\prime}(t)$. For the gaussian 
$e^{-t^2/2}$, the relevant equation is ${\frac{{\rm d}x}{{\rm d}t}} = -t e^{-t^2/2}$, 
and hence the corresponding hyperreal (fat) extension is given by ${\bf t}^2= t^2 + 
\epsilon (t) \ln t^2$. The extension of the linear variable $t$ is given by ${\bf t}= t 
+ \epsilon (t) \ln t$. The extended exponential $e^{\bf t} = t^{\epsilon} e^t$ would 
therefore have a slowly fluctuating power law tail. Note that the variable $t$ gets the 
infinitesimal correction term in the logarithmic variable $\ln t$, when the 
infinitesimal $\epsilon$ satisfies the equation 
 
$$
{\sigma{\frac{{\rm d}\epsilon}{{\rm d}\sigma}}= - \epsilon,\,\sigma = \ln t}\eqno{({\bf 
A}1)}
$$

\end{document}